\documentclass{amsart}
\usepackage{amssymb}

\newtheorem{theorem}{Theorem}[section]
\newtheorem{lemma}[theorem]{Lemma}
\newtheorem{proposition}[theorem]{Proposition}
\newtheorem{corollary}[theorem]{Corollary} 
\theoremstyle{definition}  
\newtheorem{definition}[theorem]{Definition}

\newtheorem{conjecture}[theorem]{Conjecture}  
\newtheorem{remark}[theorem]{Remark}

\newcommand{\C}{{\mathcal C}}
\newcommand{\T}{{\mathcal T}}
\newcommand{\cA}{{\mathcal A}}
\newcommand{\cI}{{\mathcal I}}

\newcommand{\cD}{{\mathcal D}}
\newcommand{\M}{{\mathcal M}}
\newcommand{\be}{{\bf 1}}
\newcommand{\bG}{{\bf \Gamma}}
\newcommand{\bV}{{\bf V}}

\newcommand{\BZ}{{\mathbb Z}}

\newcommand{\BC}{{\mathbb C}}
\newcommand{\eps}{{\varepsilon}}

\newcommand{\Ve}{\mbox{Vec}}

\newcommand{\Fun}{\mbox{Fun}}

\newcommand{\id}{\mbox{id}}
\newcommand{\Tr}{\mbox{Tr}}

\newcommand{\End}{\mbox{End}}

\title{Module categories over representations of $SL_q(2)$ in the non-semisimple case}
\author{Victor Ostrik}
\email{vostrik@darkwing.uoregon.edu}
\address{Department of Mathematics, 1222 University of Oregon, Eugene OR 97403-1222}
\thanks{The author was partially supported by NSF grant DMS-0098830.}
\date{September 2005}
\dedicatory{To Joseph Bernstein with admiration}
\begin{document}
\begin{abstract} We classify semisimple module categories over the tensor category
of representations of quantum $SL(2)$. 
\end{abstract} 
\maketitle
\section{Introduction}
Let $k$ be an algebraically closed field and let $\C$ be a tensor category over $k$. 
It is an important and interesting question to classify all semisimple module categories
over $\C$. For example, this in principle allows to construct all weak Hopf algebras
$H$ such that the category of comodules over $H$  is tensor equivalent to $\C$, that is, all 
realizations of $\C$.
Such a classification is available in some cases: when $\C$ is a group-theoretical fusion
category (for example $\C$ is the category of representations of a finite group, or a 
possibly twisted Drinfeld double of a finite group), see \cite{O2}; when $k=\BC$ and 
$\C$ is a fusion category attached to quantum $SL(2)$, see \cite{Oc,BEK,KO,O1,EO};
when $k=\BC$, $\C=\C_q$ is the category of representations of quantum $SL_q(2)$ and $q$
is not a root of unity, see \cite{EO}. The last case is the starting point for this note which
is a continuation of \cite{EO}. In this note we extend the results of \cite{EO} to the
case of arbitrary $k$ (thus allowing $k$ to have a positive characteristic) and arbitrary
$q$ (thus allowing the case when $q$ is a root of unity). In other words, we extend 
the results of \cite{EO} to the case when the tensor category in question is no longer
semisimple. It turns out that the results here are almost the same as in \cite{EO} but
the proof requires new technical tools (derived category) and explicit information on the
simple representations of quantum $SL(2)$ (tensor product theorem).

As an application of the main result of this note, we deduce, following \cite{MOV}, the Koszulity of 
preprojective algebras in some new cases. In another direction,
we give an alternative proof to some results of J.~Bichon \cite{Bi} which from our point of view
can be interpreted as a determination of the fiber functors $\C_q \to \Ve$ (see also \cite{Y2}).

I am grateful to Pavel Etingof for useful comments. Thanks are also due to the referee for careful
reading of the paper and useful comments.

\section{Main Theorem}
\subsection{Quantum $SL(2)$}\label{prel}
Let $k$ be an algebraically closed field of arbitrary characteristic. Let $q\in k^*$ be a nonzero scalar. Recall (see e.g. \cite{K}) that the Hopf algebra $SL_q(2)$ 
is defined by generators $a,b,c,d$ and relations:
$$ba=qab,\; db=qbd,\; ca=qac,\; dc=qcd,\; bc=cb,$$
$$ad-da=(q^{-1}-q)bc,\; ad-q^{-1}bc=1,$$
$$\Delta \left( \begin{array}{cc}a&b\\ c&d\end{array}\right) =
\left( \begin{array}{cc}a&b\\ c&d\end{array}\right)\otimes
\left( \begin{array}{cc}a&b\\ c&d\end{array}\right),$$
$$\eps \left( \begin{array}{cc}a&b\\ c&d\end{array}\right) =
\left( \begin{array}{cc}1&0\\ 0&1\end{array}\right),\;
S \left( \begin{array}{cc}a&b\\ c&d\end{array}\right) =
\left( \begin{array}{cc}d&-qb\\ -q^{-1}c&a\end{array}\right).$$
Let $\C_q$ denote the tensor category of finite dimensional comodules
over $SL_q(2)$ (this notation differs from the notation of \cite{EO}!)

It is well known that the category $\C_q$ is not semisimple only in the following two cases:

1) $q\ne \pm 1$ is a root of unity;

2) $q=\pm 1$ and $char (k)>0$.

In the first case we set $l$ to be the smallest positive integer such that $q^l=\pm 1$ and in the
second case we set $l=char (k)$. Recall that the Frobenius morphism defined in the two cases
above is the imbedding of Hopf algebras $SL_{q^{l^2}}(2)\subset SL_q(2)$ given by the
formulas $\bar a\mapsto a^l$, $\bar b\mapsto b^l$, $\bar c\mapsto c^l$, $\bar d\mapsto d^l$
(here $\bar a$, $\bar b$, $\bar c$, $\bar d$ are the generators of the algebra $SL_{q^{l^2}}(2)$). 
In particular we have a fully faithful tensor functor $Fr: \C_{q^{l^2}}\to \C_q$.

\subsection{Statement of the Main Theorem}
Let $\M$ be a semisimple category over $k$ with finitely many simple objects.
Recall that for any abelian tensor category $\C$ the structure of a module category over
$\C$  on $\M$ is just an exact tensor functor $F:\C \to \Fun(\M,\M)$ where $\Fun(\M,\M)$
is the tensor category of additive functors from $\M$ to itself, see \cite{O1}. If the isomorphism
classes of simple objects in $\M$ are labeled by a finite set $I$, then the category
$\Fun(\M,\M)$ is identified with the category of $I\times I-$graded vector spaces
endowed with "matrix" tensor product, see \cite{EO}. We say that $I\times I-$graded
vector space $V=\bigoplus_{i,j\in I}V_{ij}$ is {\em symmetric} if $\dim V_{ij}=\dim V_{ji}$.  
To any symmetric $I\times I-$graded vector space we attach a graph $\Gamma$ with the set
of vertices $I$ and the vertices $i$ and $j$ joined by $\dim V_{ij}$ edges.

Here is the main result of this note:

\begin{theorem}\label{main} {\em (i)} The semisimple module categories with 
finitely many simple objects over the category $\C_q$ are classified by the following data:

1) A finite set $I$;

2) A symmetric $I\times I-$graded vector space $V=\bigoplus_{i,j\in I}V_{ij}$ such that 
the corresponding graph does not contain a connected component of ADET type;

3) A collection of nondegenerate bilinear forms $E_{ij}: V_{ij}\otimes V_{ji}
\to k$,

\noindent subject to the following condition: for each $i\in I$ we have
$$\sum_{j\in I}\Tr (E_{ij}(E_{ji}^t)^{-1})=-q-q^{-1}.\eqno(*)$$
\end{theorem}

\begin{remark} It is surprising that Theorem \ref{main} has almost exactly the same formulation as 
Theorem 2.5 in \cite{EO} where it was assumed that $k=\BC$ and $q$ is not a root of unity.
Note that in {\em loc. cit.} there is no condition on the absence of the ADET type components
since in the presence of such a component equation (*) implies that $q$ is a root of unity.
\end{remark}

\begin{remark} See \cite{EO} for the list of graphs of $ADET$ type. The complex solutions of the
equation $(*)$ for these graphs also correspond to module categories but over 
the semisimple subquotient of the category $\C_q$, not over the category $\C_q$ itself, see
\cite{EO}.
\end{remark}

\subsection{Tilting modules and the universal property}

Let $\be \in \C_q$ denote the unit object and
let $\bV \in \C_q$ be a two-dimensional comodule $\bV$ with 
the basis $x,y$ and the coaction given by 
$$\Delta_\bV \left( \begin{array}{c}x\\ y\end{array}\right) =
\left( \begin{array}{cc}a&b\\ c&d\end{array}\right) \otimes
\left( \begin{array}{c}x\\ y\end{array}\right) .$$
We will call $\bV \in \C_q$ the standard object.
Recall that the object $\bV \in \C_q$ is irreducible and selfdual. Moreover for any
isomorphism $\phi :\bV \to \bV^*$ the composition
$$\be \stackrel{coev_\bV}{\longrightarrow}\bV \otimes \bV^*
\stackrel{\phi \otimes \phi^{-1}}{\longrightarrow}\bV^*\otimes \bV
\stackrel{ev_\bV}{\longrightarrow}\be \eqno(1)$$
equals to $-(q+q^{-1})\id_\be$. We fix a choice of an isomorphism $\phi :\bV \to \bV^*$
from now on.

Recall here (see e.g. \cite{D}, 1.8) that an additive category is called {\em pseudo-abelian} 
if any idempotent endomorphism is a projection on a direct factor. For any $k-$linear category
one defines its pseudo-abelian envelope; this is a category objects of which are formal
direct summands of formal finite direct sums of objects of initial category, see {\em loc. cit.} 

Let $\T_q\subset \C_q$ be the full additive subcategory of {\em tilting modules}, that is the
smallest pseudo-abelian subcategory of $\C_q$ containing the objects $\bV^{\otimes^n},\; n\in \BZ_+$
(thus $\T_q$ contains $\be =\bV^{\otimes^0}$ and $\bV=\bV^{\otimes^1}$), see \cite{An}. 
The category $\T_q$ is not  abelian in general. It is obvious that $\T_q$ is tensor subcategory of $\C_q$. 

\begin{theorem}\label{turaev1}
The triple $(\T_q,\bV,\phi)$ has the following universal property: let 
$\cD$ be a pseudo-abelian monoidal category, let $W\in \cD$ be a 
right rigid object and $\Phi: W\to W^*$ be an isomorphism such that the
composition morphism
$$\be \stackrel{coev_W}{\longrightarrow}W\otimes W^*
\stackrel{\Phi \otimes \Phi^{-1}}{\longrightarrow}W^*\otimes W
\stackrel{ev_W}{\longrightarrow}\be \eqno(2)$$
equals to $-(q+q^{-1})\id_\be$.
Then there exists a unique tensor functor $\tilde F: \T_q\to \cD$ such
that $\tilde F(\bV)=W$ and $\tilde F(\Phi )=\phi$.
\end{theorem}

\begin{proof} Let $TL(-q)$ be the Temperley-Lieb category (see \cite{Tu,B,GW,Y}, in \cite{Tu} it is
called skein category). By definition the category $TL(-q)$ has an object $X$ and two maps
$\alpha :\be \to X\otimes X$ and $\beta :X\otimes X\to \be$ such that $(\beta \otimes id_X)\circ
(id_X\otimes \alpha)=(id_X\otimes \beta)\circ (\alpha \otimes id_X)=id_X$ and $\beta \circ \alpha =
(-q-q^{-1})id_\be$. Moreover, by definition the category $TL(-q)$ is universal category with
such an object: for any tensor category $\cD$ with an object $W$ and the maps $\tilde \alpha :
\be \to W\otimes W$, $\tilde \beta :W\otimes W\to \be$ satisfying the same identities we have a
unique tensor functor $F: TL(-q)\to \cD$ such that $F(X)=W$, $F(\alpha)=\tilde \alpha$, 
$F(\beta)=\tilde \beta$. In particular we have such a functor $F: TL(-q)\to \T_q$ such that
$F(X)=\bV$, $F(\alpha)=(id_\bV \otimes \phi)\circ coev_\bV$, $F(\beta)=
ev_\bV \circ (\phi^{-1}\otimes id_\bV)$.
Now the quantum Schur-Weyl duality (we need a version over $\BZ [q,q^{-1}]$ established in
\cite{DPS}) states that the functor $F$ is fully faithful and hence the category $\T_q$ is
equivalent to the pseudo-abelian envelope of $TL(-q)$. Therefore the category $\T_q$ has the same
universal property as $TL(-q)$ but with respect to pseudo-abelian tensor categories. Finally,
the universal property involving $\alpha$ and $\beta$ is the same as the universal property
involving $\phi$: we can express $\alpha$, $\beta$ in terms of $\phi$ ($\alpha=(id_\bV \otimes \phi)\circ coev_\bV$, $\beta=ev_\bV \circ (\phi^{-1}\otimes id_\bV)$) and vice versa 
($\phi =(\beta \otimes id_{\bV^*})\circ (id_\bV \otimes coev_\bV)$ and $\phi^{-1}=
(id_\bV \otimes \alpha^*)\circ (coev_\bV \otimes id_{\bV^*})$). The Theorem is proved.

\end{proof} 

\begin{corollary}\label{tilting}
 Let $\M$ be a semisimple category with isomorphism classes of simple
objects labeled by a finite set $I$. The tensor functors $\tilde F: \T_q\to \Fun(\M,\M)$ are
classified by the following data: 

1) A symmetric $I\times I-$graded vector space $V=\bigoplus_{i,j\in I}V_{ij}$;

2) A collection of nondegenerate bilinear forms $E_{ij}: V_{ij}\otimes V_{ji} \to k$,
such that the equation $(*)$ holds.
\end{corollary}

Now we explain why the graphs of ADET type are special from the point of view
of this Corollary: 

\begin{lemma} \label{ade}
Assume that the graph $\Gamma$ attached to $V$ is of ADET type.
Then there exists a polynomial $P$ such that 

1) $P([\bV])\in K(\C_q)$ represents an actual (not just virtual) non-zero object;

2) $P([\tilde F(\bV)])=0\in K(\Fun(\M,\M))$.
\end{lemma}

\begin{proof} Let $h$ be the Coxeter number of $\Gamma$. Consider the polynomial
$P(x)=P_{h-1}(x)=U_{h-1}(\frac{x}2)$ where $U_n$ is the Tchebysheff polynomial of
second kind (thus $P_n(2\cos x)=\frac{\sin (n+1)x}{\sin x}$). It is well known (see \cite{AW}) that
$P_{h-1}([\bV])=[V_{h-1}]$ where $V_{h-1}$ is the Weyl module  and thus
1) holds. On the other hand, it is well known that all the eigenvalues of the matrix
$\dim (V_{ij})$ (= adjacency matrix of $\Gamma$) are of the form $2\cos \frac{\pi m_i}{h}$
for some integers $m_i$, $0<m_i<h$; in particular all the eigenvalues are
roots of $P_{h-1}$.  Since the matrix $\dim (V_{ij})$ is diagonalizable we get that
$P([\tilde F(\bV)])=P(\dim (V_{ij}))=0\in K(\Fun(\M,\M))$ and 2) holds. Lemma is proved.
\end{proof}

Here is an outline of the proof of Theorem \ref{main}: for a tensor functor $F: \C_q \to \Fun(\M,\M)$
by restriction we get a functor $\tilde F:  \T_q\to \Fun(\M,\M)$ and hence by Corollary \ref{tilting}
a symmetric $I\times I-$graded vector space and a solution of the equation $(*)$. It follows
from Lemma \ref{ade} that the graph $\Gamma$ attached to $V$ does not contain a connected
component of ADET type (since a tensor functor can not send a nonzero object to zero). The real
difficulty is in the proof of Theorem \ref{main} in converse direction: let $(V, E_{ij})$ be a solution
to the equation $(*)$ such that the graph $\Gamma$ does not contain a connected component of
ADET type; we need to show that the corresponding functor $\tilde F: \T_q\to \Fun(\M,\M)$ extends 
uniquely to a tensor functor $F: \C_q\to \Fun(\M,\M)$.  This will be done in subsequent sections.

\subsection{Derived category} In this section we show that if the extension of the functor
$\tilde F: \T_q\to \Fun(\M,\M)$ exists then it is unique. For this let us consider two
triangulated categories: the homotopy category $K^b(\T_q)$ of bounded complexes in $\T_q$ 
and the bounded derived category $D^b(\C_q)$. The inclusion $\T_q\subset \C_q$ induces
an exact functor $A: K^b(\T_q)\to D^b(\C_q)$. Note that $A$ has an obvious structure of tensor
functor. The following observation is crucial for this paper:

\begin{proposition} (\cite{BBM}) The functor $A$ is an equivalence of categories.
\end{proposition}

\begin{proof} It is known that $\C_q$ is a highest weight category (in a sense of \cite{CPS}), 
see \cite{AW}. For every highest weight category one defines the class of tilting modules (see
\cite{An} for definition in the case of quantum groups, which applies to $\C_q$). It is known from
\cite{An} that for the category $\C_q$ the definition of tilting modules from Section 2.3 agrees
with this general definition for highest weight categories. Finally it is explained in \cite{BBM} 1.5
that for any highest weight category $\C$ the bounded homotopy category of tilting modules 
is equivalent to the bounded derived category $D^b(\C)$.
\end{proof}

\begin{corollary} If the functor $\tilde F: \T_q\to \Fun(\M,\M)$ extends to a functor $F: \C_q\to \Fun(\M,\M)$
then an extension is unique up to an isomorphism of tensor functors.
\end{corollary}

\begin{proof} The functor $\tilde F$ extends uniquely to the exact functor $R\tilde F: K^b(\T_q)\to
 K^b(\Fun(\M,\M))=D^b(\Fun(\M,\M))$. Let $B: \C_q\to D^b(\C_q)$ be the obvious imbedding and
let $H^i: D^b(\Fun(\M,\M))\to \Fun(\M,\M)$ be the $i-$th cohomology functor. It is clear
that the functor $\tilde F$ extends if and only if $H^i\circ R\tilde F\circ A^{-1}\circ B=0$ for 
$i\ne 0$ and the extension is isomorphic to $H^0\circ R\tilde F\circ A^{-1}\circ B$. 
 \end{proof}

\subsection{Key Lemma}\label{key}
 Our goal now is to prove the vanishing $H^i\circ R\tilde F\circ A^{-1}\circ B=0$
as above. For this we need some explicit information on the simple objects in the category $\C_q$.
The only nontrivial case is when the category $\C_q$ is not semisimple.
Thus we assume that 

1) either $q\ne \pm 1$ is a root of unity and $l$ is the smallest positive integer such that
$q^l=\pm 1$;

 2) or $q=\pm 1$ and $k$ has characteristic $l>1$.
 
 The following facts are well known (and easy to prove):
 
 (a) (\cite{AW}) The simple objects in $\C_q$ are labeled by its highest weight which 
 is an arbitrary nonnegative integer. We will denote the simple object with highest 
 weight $k\in \BZ_+$ by $L_k$. In particular $L_0=\be$, $L_1=\bV$.
 
 (b) (\cite{An}) The simple objects $L_0, L_1, \ldots, L_{l-1}$ are tilting. Moreover, we have
 an equality in the Grothendieck ring $[L_m]=P_m([\bV])\in K(\C_q)\; ,m=0,\ldots, l-1 $ where $P_m$
 is defined by $P_m(2\cos x)=\frac{\sin (m+1)x}{\sin x}$ (thus $P_m(2x)$ is the Tchebysheff
 polynomial of second kind).
 
 (c) (\cite{An}) The object $\bV \otimes L_{l-1}$ is not semisimple; it has length 3; its socle and 
 cosocle are both isomorphic to $L_{l-2}$ and the third simple constituent is isomorphic to 
 $L_l$. Let $f: L_{l-2}\to \bV\otimes L_{l-1}$ and $g: \bV\otimes L_{l-1}\to L_{l-2}$ be the 
 corresponding maps; we may and will assume that $g=f^*$. This also implies that in
 the Grothendieck ring $[L_l]=Q_l([\bV])\in K(\C_q)$ where $Q_l$ is defined by $Q_l(2\cos(x))=
 2\cos(lx)$.  
  
  (d) (\cite{AW}) Tensor product theorem: for any integers $k\ge 1$ and $0\le m\le l-1$ the simple module 
  $L_{lk+m}$ is isomorphic to $Fr(\bar L_k)\otimes L_m$. Here $Fr$ is the Frobenius functor, see
  \ref{prel} and $\bar L_k$ is a simple object of $\C_{q^{l^2}}$ with highest weight $k$.
  
  (e) The special cases of (d): $L_l=Fr(\bar \bV)$ (here $\bar \bV$ is the standard object 
  of $\C_{q^{l^2}}$),
  and $L_{2l-1}=Fr(\bar \bV)\otimes L_{l-1}$. This also implies that $L_{2l-1}$ is tilting, see \cite{An}.
  
 \begin{lemma} \label{klemma} Assume that the graph $\Gamma$ is connected and not of ADET type. 
 We have $H^i\circ R\tilde F\circ A^{-1}\circ B(L_l)=0$ for $i\ne 0$.
 \end{lemma}
 
 \begin{proof} It follows from (c) above that $A^{-1}\circ B(L_l)$ can be represented by the
 following complex in $K^b(\T_q)$:
 
 $$\cdots \to 0\to L_{l-2}\stackrel{f}{\longrightarrow} \bV\otimes L_{l-1}\stackrel{g}{\longrightarrow} 
 L_{l-2}\to 0\to \cdots$$ 

Thus the statement of the Lemma is obvious for $i\ne \pm 1$.
This complex is obviously self-dual, so we only need to prove that 
$H^{-1}\circ R\tilde F\circ A^{-1}\circ B(L_l)=0$ or, equivalently, that the
map $F(f): F(L_{l-2})\to F(\bV \otimes L_{l-1})$ is injective. We will need the
following

{\bf Sublemma 1.} The injective map $f\otimes id: L_{l-2}\otimes L_{l-1}\to 
\bV \otimes L_{l-1}\otimes L_{l-1}$ is an imbedding of a direct summand.

\begin{proof} This is easy when $char(k)=0$ since $L_{l-1}$ is injective in this case. In general
note that the cohomology of the complex

$$\cdots \to 0\to L_{l-2}\otimes L_{l-1}\stackrel{f\otimes id}{\longrightarrow} 
\bV \otimes L_{l-1}\otimes L_{l-1}\stackrel{g\otimes id}{\longrightarrow} 
L_{l-2}\otimes L_{l-1}\to 0\to \cdots$$ 
is $L_l\otimes L_{l-1}=L_{2l-1}$ (see (e)). The modules $L_{2l-1}$ and $L_{l-2}\otimes L_{l-1}$ are both
tilting and hence $Ext^1(L_{2l-1},L_{l-2}\otimes L_{l-1})=0$ (since all higher $Ext$ groups between
tilting modules vanish, see \cite{An}). The Sublemma follows.
\end{proof}

Sublemma 1 implies that the map $F(f)\otimes F(id): F(L_{l-2})\otimes F(L_{l-1})\to 
F(\bV \otimes L_{l-1})\otimes F(L_{l-1})$ is injective. Thus if $M=\mbox{Ker} F(f)$ then 
$M\otimes F(L_{l-1})=0$.

{\bf Sublemma 2.} If $\Gamma$ is connected and not of ADET type then $id$ is a direct summand of 
$F(L_{l-1})\otimes F(L_{l-1})$ in $\Fun(\M,\M)$.

\begin{proof} Let $\M_\BC$ be a semisimple category over $\BC$ such that the Grothendieck groups
$K(\M)$ and $K(\M_\BC)$ are isomorphic. We identify the Grothendieck rings $K(\Fun(\M,\M))$ and
$K(\Fun(\M_\BC,\M_\BC))$ using such an isomorphism. It is known from \cite{EO} that under our
assumptions  there is $q_\BC \in \BC^*$  which is not a root of unity and the tensor functor $F_\BC: \C_{q_\BC}\to \Fun(\M_\BC,\M_\BC)$ such that $[F_\BC(\bV^\BC)]=[F(\bV)]\in K(\Fun(\M,\M))$ (here 
$\C_{q_\BC}$ is considered over $\BC$, not over $k$, and $\bV^\BC$ is the standard object of 
$\C_{q_\BC}$). This implies that $[F_\BC(L_{l-1}^\BC)]=[F(L_{l-1})]$ (here $L_{l-1}^\BC$ is the object of $\C_{q_\BC}$ with highest weight $l-1$). But it is known that $\be^\BC$ (the unit object of 
$\C_{q_\BC}$) is a direct summand of $L_{l-1}^\BC \otimes L_{l-1}^\BC$. The Sublemma follows.
\end{proof}
 
Sublemma 2 implies that $M$ is a direct summand of $M\otimes F(L_{l-1})\otimes F(L_{l-1})=0$
and hence that $M=0$. The Lemma is proved.

\end{proof}

\subsection{Proof of the Main Theorem} We just need to prove that under the assumptions of
the Theorem \ref{main} the functor $\tilde F:\T_q \to \Fun(\M,\M)$ extends to the functor
$F: \C_q \to \Fun(\M,\M)$ or, equivalently, that  $H^i\circ R\tilde F\circ A^{-1}\circ B(N)=0$
for $i\ne 0$ and any $N\in \C_q$. We can restrict ourselves to the case when
$\Gamma$ is connected (otherwise $F$ is just a direct sum of functors corresponding to the
connected components of $\Gamma$). 

First we prove that $H^i\circ R\tilde F\circ A^{-1}\circ B(L_m)=0$
for $i\ne 0$ and simple $L_m$. By \ref{key} (b) we know that this is true for $m=0,\ldots ,l-1$. Thus
by \ref{key} (d) we are reduced to the simple modules $Fr(\bar L_k)$. By Lemma \ref{klemma}
the desired vanishing is known for $Fr(L_1)=L_l$. Moreover, the object 
$W=H^0\circ R\tilde F\circ A^{-1}\circ i(L_l)\in \Fun(\M,\M)$ is endowed with isomorphism 
$\psi: W\to W^*$ such that the composition $\be \to W\otimes W^*\stackrel{\psi \otimes 
\psi^{-1}}{\longrightarrow}W^*\otimes W\to \be$ equals to $-(q^{l^2}+q^{-l^2})=\pm 2$ (this comes
from any choice of an isomorphism $L_l\to L_l^*$). Thus we can apply the same machinery to $W$
as we applied before (see Corollary \ref{tilting}) to $V$. This completes the proof if $char(k)=0$
since in this case the category $\C_{q^{l^2}}$ is semisimple. In general, we are going to apply
the induction and for this we need to know that the graph $\bar \Gamma$ associated with $W$
has no connected components of ADET type.  Assume it does. Then there exists a polynomial 
$P_n$ (see Lemma 2.6) and simple object $M\in \M$ such that $P_n([W])M=0\in K(\M)$. On the other
hand it is known that in $K(\C_q)$ we have an equality $[L_{l-1}]P_n([L_l])=P_{ln+l-1}([V])$
(equivalently, $\frac{\sin(lx)}{\sin(x)}\frac{\sin((n+1)lx)}{\sin(lx)}=\frac{\sin((n+1)lx)}{\sin(x)}$).
Thus passing to the category $\C_{q_\BC}$ and $\M_\BC$ (see the proof of Sublemma 2) 
we will have the objects $L^\BC =L_{ln+l-1}^\BC \in \C_{q_\BC}$ and $M^\BC \in \M_\BC$ such
that $F_\BC (L^\BC)M^\BC =0$. But this is impossible since $L^\BC \otimes L^\BC$ containd
$\be^\BC$ as a direct summand. Thus it is proved by induction in the highest weight $m$ that
$H^i\circ R\tilde F\circ A^{-1}\circ B(L_m)=0$ for $i\ne 0$ and simple $L_m\in \C_q$.

The long exact sequence in cohomology implies now that the set of $N\in \C_q$ such that
$H^i\circ R\tilde F\circ A^{-1}\circ B(N)=0$ for $i\ne 0$ is closed under extensions and hence
coincides with the set of all objects in $\C_q$. Thus the functor $F(N)=H^0\circ R\tilde F\circ A^{-1}\circ B(N)$ is the desired extension of the functor $\tilde F$ (note that the functor $F$ has a structure
of tensor functor since it is a composition of tensor functors). Theorem \ref{main} is proved.

\begin{remark}  Assume that $k=\BC$.
A similar argument proves that for any abelian tensor category $\cA$ and a tensor functor
$\tilde F: \T_q \to \cA$ such that $F(L_{l-1})\ne 0$  we have an extension $F: \C_q\to \cA$
(we assume that $\End(\be_\cA)=k$). Recall that $L_{l-1}$ generates the tensor ideal 
$\cI_q\subset \T_q$ of negligible modules, see \cite{An} and the quotient $\tilde \C_q=\T_q/\cI_q$ is
semisimple. Thus any tensor functor $F:\T_q \to \cA$ either extends
to $\C_q$ or factors through the semisimple quotient $\tilde \C_q$. We can say that the category
$\T_q$ admits exactly two {\em abelian extensions}.
It would be interesting to investigate
similar questions for other tensor categories defined by universal properties, in 
particular for Deligne's categories from \cite{D}.
\end{remark}

\section{Some applications}
\subsection{Preprojective algebras} Recall (see \cite{DR}) that a {\em modulated graph} $\bG$ is 
a finite set $I$, collection of finite dimensional vector spaces $\{ V_{ij}\}_{i,j\in I}$ and a collection of nondegenerate bilinear forms $E_{ij}: V_{ij}\otimes V_{ji}\to k$ (the definition in \cite{DR} is
slightly more general); the underlying graph $\Gamma$ of a modulated graph $\bG$ has the set
of vertices $I$ and $\dim(V_{ij})$ edges joining $i$ and $j$. In other words we have $\bG$ is the 
same as $I\times I-$graded vector space $V$ and a nondegenerate pairing $E: V\otimes V\to \be$. 
Let $E^*: \be \to V\otimes V$ be the dual of $E$. Let $T^\bullet (V)$ be a tensor algebra of $V$;
that is $T^n(V)=V^{\otimes^n}$ with obvious multiplication. The algebra $T(V)$ identifies with
{\em path algebra} of $\Gamma$.

\begin{definition} (\cite{DR})  The preprojective algebra $P(\bG)$ associated with $\bG$ is the
quotient $T(V)$ by the ideal generated by the image of $E^*:\be \to T^2(V)$.
\end{definition}

Note that by definition $P(\bG)$ is graded.
We don't know a counterexample to the following 

\begin{conjecture}\label{con} (i) Let us fix a connected graph $\Gamma$ with the adjacency matrix $A$. 
The family of preprojective algebras $P(\bG)$ is flat; that is for any choice of $k$ and $E$ 
the matrix Hilbert series (see \cite{MOV}) is constant. In particular, if $\Gamma$ is
not of ADET type the Hilbert series equals to $(1-At+t^2)^{-1}$.

(ii) Assume that $\Gamma$ is not of ADET type. Then the algebras $P(\bG)$ are Koszul.
\end{conjecture}

Conjecture \ref{con} is known to hold in a number of cases. For example it is known
that part (i) holds for the graphs of ADE type (this follows from an interpretation of preprojective
algebras in terms of representations of quivers due to Gelfand-Ponomarev, see e.g. \cite{DR})
and it is probably easy for graphs of T type.
Theorem \ref{main} implies Conjecture \ref{con} in some cases. Recall that for a matrix
$(b_{ij})$ the eigenvalue $\lambda$ is called nondegenerate if there exists a 
$\lambda-$eigenvector $(r_i)_{i\in I}$  such that $\prod_{i\in I}r_i\ne 0$. To a
modulated graph $\bG$ we associate the matrix $D=(d_{ij})$ where $d_{ij}=\Tr(E_{ij}(E_{ji}^t)^{-1})$.

\begin{proposition} Assume that $\Gamma$ is connected and not of ADET type.
Assume that the matrix $D$ admits a nondegenerate eigenvalue $\lambda$. Then the algebra
$P(\bG)$ is Koszul and its Hilbert series is given in Conjecture \ref{con}.
\end{proposition}

\begin{proof} Let $r$ be an automorphism of $\be$ (in other words $r$ is just a collection
of nonzero scalars $r_i$). Consider new modulated graph $\bG'$ obtained from $\bG$ by
replacing $E$ by $E': \be \stackrel{r}{\longrightarrow} \be \stackrel{E}{\longrightarrow} V\otimes V$.
It is easy to see that the algebras $P(\bG)$ and $P(\bG')$ are canonically isomorphic (they have the
same generators and relations differ just by scalars). Obviously $E'_{ij}=r_jE_{ij}$. Hence
for the matrix $D'$ attached to $E'$ we have $d_{ij}'=r_jd_{ij}r_i^{-1}$. Thus $\sum_jd'_{ij}=\lambda$.
Hence we can construct a tensor functor $F: \C_q\to \Fun(\M,\M)$ using $E'$ and Theorem \ref{main}
(here $q$ is determined from the equation $q+q^{-1}=-\lambda$).
The rest of the proof is the same as in \cite{MOV}: we identify $P(\bG')$ with $F(S_q)$, where
$S_q\in \C_q$ is the $q-$symmetric algebra.
\end{proof}

\begin{remark} Unfortunately the conditions of this Proposition are too restrictive; for example
it gives no answer in the case when $\Gamma$ is a star-shaped tree with $n+1$ vertices
and $char(k)$ divides $n$.
\end{remark}

\begin{remark} New results in the direction of Conjecture \ref{con} are contained in a recent
preprint \cite{EE}. The authors of {\em loc. cit.} proved Conjecture \ref{con} in many new
cases using a completely different elementary method. Moreover it seems likely that the
methods of {\em loc. cit.} would imply Conjecture \ref{con} in complete generality.
\end{remark}

\subsection{Some Hopf algebras} The following class of Hopf algebras was defined in
\cite{DL}. Let $V$ be a finite dimensional vector space endowed with a nondegenerate 
bilinear form $E$. Let $e_\alpha$ be a basis of $V$. Consider the algebra $H(E)$
generated by the elements $a_{\alpha \beta}$ with the following relations given in
matrix form:

$$E^{-1}a^tEa=aE^{-1}a^tE=I$$
where $I$ is the identity matrix. The algebra $H(E)$ has the following structure of Hopf algebra:

$$\delta (a)=a\otimes a, \eps (a)=I, S(a)=E^{-1}a^tE.$$

Observe that the space $V$ has an obvious structure of comodule over $H(E)$:
$\delta_V(e_\alpha)=\sum_\beta e_\beta \otimes a_{\beta \alpha}$. Moreover the
map $E: V\otimes V\to k$ is a comodule map. Actually the Hopf algebra $H(E)$ is universal
with these properties, see \cite{Bi} Proposition 2.2 (ii).

The corepresentation theory of the algebras $H(E)$ was determined by J.~Bichon \cite{Bi}.
Here we reprove his result from our point of view.

\begin{theorem} (\cite{Bi}) The tensor category of finite dimensional comodules over $H(E)$
is tensor equivalent to $\C_q$ where $q$ is determined from the equation $\Tr (E(E^t)^{-1})=-q-q^{-1}$.
\end{theorem}

\begin{proof} Theorem \ref{main} with $|I|=1$ states that the tensor functors $K: \C_q \to \Ve$
(here $\Ve$ is the tensor category of vector spaces) are classified by the pairs $(V,E)$ as above
with $\Tr (E(E^t)^{-1})=-q-q^{-1}$. 

Now the universal property of $H(E)$ can be restated as the following universal property of
the category $H(E)-comod$ of $H(E)-$comodules: for any abelian tensor category $\cD$ with 
the fiber functor
$G: \cD \to \Ve$, the object $W\in \cD$ endowed with bilinear pairing $E_1 : W\otimes W\to \be$
such that $G(W)=V$ and $G(E_1)=E$ there exists a tensor functor $F: H(E)-comod \to \cD$ such
that $F(V)=W$, $F(E)=E_1$ and the composition $G\circ F$ is isomorphic to the forgetful functor
$H(E)-comod \to \Ve$. 

Finally using the universal property of the category $\C_q$ (or rather $\T_q$) we can construct
a tensor functor $\C_q\to H(E)-comod$ and using the universal property of the category
$H(E)-comod$ and the functor $K$ we can construct a tensor functor $H(E)-comod \to \C_q$.
It is obvious that these functors are mutually inverse. The Theorem is proved.

\end{proof}

\begin{remark} Using the module categories constructed in Theorem \ref{main} we can 
similarly construct some weak Hopf algebras with corepresentation categories $\C_q$.
\end{remark}

\end{document}